\documentclass[12pt]{amsart}
\advance\hoffset by -0.5 truecm

\usepackage{amsmath}
\usepackage{amssymb}
\usepackage{amsthm}

\newtheorem{theorem}{Theorem}

\newtheorem{proposition}{Proposition}

\newtheorem{corollary}{Corollary}
\begin{document}

\title{Ricci-positive geodesic flows and point-completion of static monopole fields}
\author[K. Dorji, A. Harris]{Kumbu Dorji and Adam Harris}
\address{School of Science and Technology, University of New England,
Armidale, NSW 2351, Australia}
\email{adamh@une.edu.au\\kdorji3@myune.edu.au}

\begin{abstract} 
Let $(\hat{M}, g)$ be a compact, oriented Riemannian three-manifold corresponding to the metric point-completion $M\cup\{P_{0}\}$ of a manifold $M$, and let $\xi$ denote a geodesible Killing unit vector field on $\hat{M}$ such that the Ricci curvature function $Ric_{g}(\xi) > 0$ everywhere, and is constant outside a compact subset $K\subset\subset M$. Suppose further that $(E,\nabla, \varphi)$ supply the essential data of a monopole field on $M$, smooth outside isolated singularities all contained in $K$. The main theorem of this article provides a sufficient condition for smooth extension of $(E, \nabla, \varphi)$ across $P_{0}$, in terms of the {\em Higgs potential} $\Phi$, defined in a punctured neighbourhood of $P_{0}$ by  
\[\nabla_{\xi}\Phi - 2{\bf i}[\varphi, \Phi] = \varphi \ .\]
The sufficiency condition is expressed by a system of equations on the same neighbourhood, which can be effectively simplified in the case that $\hat{M}$ is a regular Sasaki manifold, such as the round 
${\Bbb S}^{3}$. 

\end{abstract}

\date{June 2018}
\maketitle

\section{Introduction}

In its formal conception, every monopole on a three-dimensional manifold $M$ is comprised of a smooth complex vector bundle $E$ over $M$ (or a subset of $M$, if singularities are present) with structure group ${\Bbb S}{\Bbb U}(n)$, a connection $\nabla$, and an endomorphism $\varphi$ (the {\em Higgs field}) with coefficients in the corresponding Lie Algebra. A Riemannian metric $g$ and associated volume form are also needed in order to define a Hodge-star operator on differential forms. These structures are then related via the Bogomolny equation
\[ 2\star\nabla\varphi = F_{\nabla} \  \ ,\]
where $F_{\nabla}$ denotes the curvature two-form. On one hand, if the metric completion $\hat{M}$ of $M$ is compact, it is known that a non-trivial, non-singular solution of this equation cannot exist on $\hat{M}$. On the other, if $\hat{M}$ is not compact, non-singular solutions may exist provided certain asymptotic conditions are met at infinity. When $\hat{M}$ corresponds to Euclidean or Hyperbolic three-space, on which the bundle structure $E$ is effectively trivial if singularities are absent, these conditions were precisely formulated as part of a theory of moduli of solutions of the Bogomolny equation initially developed by Atiyah \cite{A3}, Hitchin \cite{Hi}, Ward \cite{W} and their collaborators. 

By contrast, the concern of the present note will be with monopoles over three-manifolds $M$ having a single-point compact Riemannian metric completion $\hat{M} = M\cup\{P_{0}\}$. It will be assumed that the monopole field has singularities at isolated points $p_{i}$ all contained within a compact subset $K\subset\subset M$. Moreover, following \cite{BiHu} and much of the tradition concerning singular monopoles, these may be assumed to be of {\em Dirac type}, which entails that in a neighbourhood $B_{\varepsilon}(p_{i})$ of each singularity, the bundle $E\mid_{
B_{\varepsilon}(p_{i})\setminus\{p_{i}\}}$ splits as a direct sum of non-trivial complex line bundles. Setting aside their precise nature, the specific problem addressed here will be to formulate a sufficient condition for smooth extension of $\varphi$ and $\nabla$ from $M\setminus K$ to 
$\hat{M}\setminus K$, thereby making a {\em removable singularity} of $P_{0}$. In particular, such a condition must first imply the existence of a smooth trivialization of the bundle $E$ in a punctured neighbourhood of $P_{0}$. Several years ago this problem was tentatively addressed by one of the present authors \cite{H} in the context of Euclidean monopoles, but our renewed interest has been inspired by recent work of Biswas and Hurtubise \cite{BiHu}, who have extended the investigation of monopoles to Sasakian three-folds. It is well-known that an anti self-dual structure on the pullback $\pi^{*}E$ to the product $M\times(0,\infty)$ is naturally supported by a Sasakian structure on $\hat{M}$. As a result $\pi^{*}E$ inherits a fully integrable holomorphic structure, while $\hat{M}\times(0,\infty)$ becomes a K\"{a}hler surface $X$. A central theme of \cite{BiHu} is to develop the machinery of the Kobayashi-Hitchin correspondence for monopoles in the regular Sasakian setting by first compactifying the fibres of $X$ and working instead with the properties of an available Gauduchon metric. The goal of the present article, however, is simply to use local tools of complex analysis on the K\"{a}hler surface $X$ in order to study removable three-dimensional point singularities. In fact, for present purposes, the strict Sasakian property is only required in a neighbourhood of $P_{0}$. 

In \cite{HP} it was shown that if $(\hat{M}, g)$ is a compact, oriented Riemannian three-manifold, with a geodesible unit vector field $\xi$ along which the Ricci tensor is strictly positive, i.e., $Ric(\xi) > 0$, then the dual one-form $\vartheta$ corresponding to the contraction $\iota_{\xi}g$ defines a contact structure on $\hat{M}$, with $\xi$ as its Reeb vector field. If $\xi^{\perp}$ denotes the distribution corresponding to $ker(\vartheta)$, then a natural endomorphism $j$ of $\xi^{\perp}$ is defined, such that for any $v\in\xi^{\perp}_{p}$ the set $\{v, jv, \xi\}$ forms a positively oriented orthonormal frame of $T_{p}\hat{M}$. As will be discussed further in the next section (cf. Proposition 1), $X$ must of course take on the structure of a Riemannian cone with respect to the metric tensor
\[\bar{g} = t^{2}g + \frac{1}{4f}dt\otimes dt \ , \ \text{for} \ t\in(0,\infty) \ , \ f:= Ric(\xi) 
\ ,\]
such that the Sasakian property is contingent on the identification 
\[d\vartheta = 2\sqrt{f} \ \iota_{j}g \ .\]
This entails two basic hypotheses : (i) that the Reeb vector field $\xi$ is {\em Killing}, hence its flow preserves $j$ as well as the volume form, as discussed in \cite{HP},  and (ii) 
that the coefficient determined by $f$ is a constant $c$, at least in a neighbourhood of $P_{0}$, where the K\"{a}hler structure is needed. It is worth mentioning at this point that if one were instead to consider a {\em Lorentzian} metric
\[  \bar{g} = t^{2}g - \frac{1}{4f}dt\otimes dt \ ,\]
as in recent work of Aazami \cite{Az}, where the Newman-Penrose formalism has been applied with interesting effect to methods of symplectic geometry, then the K\"{a}hler structure may be recovered from a simple reversal of orientation on $\hat{M}$, effectively replacing $j$ by $-j$. 

The presence of a K\"{a}hler structure on $\pi^{-1}(\hat{M}\setminus K)$ enables the standard presentation of the Bogomolny equation as a time-independent reduction of the anti self-dual property of the connection 
\[\nabla ' := d + \pi^{*}A_{\nabla} + \frac{1}{t\sqrt{c}}\pi^{*}\varphi\otimes dt \]
defined on $\pi^{*}E$, and consequently induces a holomorphic structure on $\pi^{*}E$. The first requirement for a removable singularity at $P_{0}$ is to obtain a trivialization of $E\mid_{B_{\varepsilon}(P_{0})\setminus\{P_{0}\}}$. In earlier work \cite{H} this was shown to follow from the existence of a relative holomorphic connection on $\pi^{*}E$, by means of a Hartogs technique for holomorphically extending sections across $\pi^{-1}(P_{0})$. Existence of a relative holomorphic connection is itself obstructed by complex-analytic cohomology \cite{A2}, expressible here in terms of solvability of a Cauchy-Riemann equation of the form 
\[ \bar{\partial}\Psi = \iota_{Z} F_{\nabla '}  \ ,\]
where $Z$ denotes a holomorphic vector field transverse to $\pi^{-1}(P_{0})$ (cf. section four below). Given the three-dimensional setting of our problem, however, it will be appropriate to formulate a sufficient condition for solution of this equation in terms of a $t$-independent reduction. In particular, a solution $\Psi$ should be of the form $\pi^{*}\psi$ (cf. Proposition 2). Existence of $\psi$ can be made contingent on $\varphi$, by first defining the {\em Higgs potential} $\Phi$
such that 
\[ \nabla_{\xi}\Phi - 2{\bf i}[\varphi, \Phi] = \varphi \ \]
in a punctured neighbourhood of $P_{0}$, and then setting $\psi := -2{\bf i}\nabla_{Z}\Phi$ (cf. section three). A sufficient condition for solution of the Cauchy-Riemann equation above can then be formulated in terms of the equations
\[\nabla_{\bar{Z}}(\nabla_{Z}\Phi) + \nabla_{\xi}\varphi = 0 \ ,\]
\[\nabla_{[Z,\xi]}\Phi - {\bf i}[\nabla_{Z}\varphi, \Phi] = 0 \ ,\]
(Proposition 3). We note in particular that if $(\hat{M}, g)$ corresponds to a {\em regular} Sasaki manifold, in which the Reeb flow endows $\hat{M}$ with the structure of an ${\Bbb S}^{1}$-principal bundle over a compact Riemann surface, then the Lie bracket $[Z,\xi] = 0$ (cf. \cite{BiHu}), and the equations of Proposition 3 are correspondingly simplified. It remains in section four to carry out the Hartogs extension of $\pi^{*}E$ and thereby extend $E$ smoothly across $P_{0}$, before proceeding in section five to consider extension of $\varphi$ and $\nabla$.   

{\bf Theorem}: (cf. Section five) {\em Let} $(\hat{M}, g)$ {\em be a compact, oriented Riemannian three-manifold corresponding to the point-completion of a manifold} $M$, {\em and let} $\xi$ {\em denote a geodesible Killing unit vector field on} $\hat{M}$ {\em such that the Ricci curvature function} $Ric_{g}(\xi) > 0$ {\em everywhere, and is constant outside a compact subset} $K\subset M$. {\em Suppose further that} $(E,\nabla, \varphi)$ {\em supply the essential data of a monopole field on} $M$, {\em smooth outside isolated singularities all contained in} $K$, {\em and that the following hold in the complement of} $K$:

\vspace{.1in}

{\em (i) the associated Higgs potential} $\Phi$ {\em satisfies the equations}
\[ \nabla_{\bar{Z}}(\nabla_{Z}\Phi) + \nabla_{\xi}\varphi = 0 \ ,\]
{\em and}
\[ \lambda\nabla_{Z}\Phi + [\nabla_{Z}\varphi, \Phi]  = 0\]
{\em where} $Z = \sigma - {\bf i}j\sigma$ {\em with respect to an orthonormal framing} $\{\sigma, j\sigma, \xi\}\in C^{\infty}(\hat{M}\setminus K, T\hat{M})$, {\em and a real-valued function} $\lambda = \frac{1}{4}g([\sigma, \xi], j\sigma) $,

\vspace{.1in}

{\em (ii) the Higgs field} $\varphi$ {\em and monopole connection} $\nabla$ {\em are both uniformly} $C^{1,\alpha}$ {\em with respect to admissible} $C^{1,\alpha}$-{\em framings of} $E$.

\vspace{.1in}

{\em Then the bundle} $E$ {\em extends smoothly across} $P_{0}$, {\em and admits} $C^{1}$-{\em extensions of both} $\varphi$ {\em and} $\nabla$, {\em together with a continuous extension of} $F_{\nabla}$.

\vspace{.1in}

{\bf Corollary}: {\em When} $\hat{M} = {\Bbb S}^{3}$, {\em equipped with the round metric, and the flow of} $\xi$ {\em induces the Hopf fibration, then the simplified condition}

\[(i)\qquad  \nabla_{\bar{Z}}(\nabla_{Z}\Phi) + \nabla_{\xi}\varphi
 = [\nabla_{Z}\varphi, \Phi]  = 0 \ ,\]
  
{\em together with condition} $(ii)$, {\em is sufficient for extension of} $E, \varphi, \nabla$ {\em and} $F_{\nabla}$ {\em as above.}

\vspace{.1in}

\section{Monopole fields and three-dimensional contact geometry} 

\vspace{.1in}

Let $\hat{M}$ be a smooth, orientable three-manifold corresponding to the point compactification 
of a manifold $M$, i.e., $\hat{M} = M\cup\{P_{0}\}$. We will assume that $\hat{M}$ is equipped with a 
Riemannian metric $g$ and a geodesible unit vector field $\xi\in C^{\infty}(\hat{M}, T\hat{M})$. 
Let $\vartheta\in C^{\infty}(\hat{M}, T^{*}\hat{M})$ be the metric-dual of $\xi$, and $D$ the Levi-Civita
connection associated with $g$. We note that the geodesibility requirement $D_{\xi}\xi = 0$ is 
easily seen to be equivalent to the condition $\iota_{\xi}d\vartheta = 0$. $\vartheta$ will then correspond 
to a contact form on $\hat{M}$ if $d\vartheta\mid_{\xi^{\perp}}$ is non-degenerate. A criterion of 
non-degeneracy was shown in \cite{HP}, Lemma 1, to follow from the assumption $Ric(\xi) > 0$ on $\hat{M}$, where
the function $Ric: S\hat{M}\rightarrow{\Bbb R}$ is naturally induced on the unit-sphere bundle by the 
Ricci tensor. Under this assumption we may then identify $\xi$ with the Reeb vector field of $\vartheta$.
In addition, for any $(p,\it{v})\in S\hat{M}$, there is a unique linear map $j_{(p,\it{v})}:\it{v}^{\perp}\rightarrow\it{v}^{\perp}\subset T\hat{M}$ such that $j_{(p,\it{v})}^{2} = -id\mid_{\it{v}^{\perp}}$, and 
for any unit vector $\it{u}\in\it{v}^{\perp}$, the triple $\{u, ju, v\}$ is a positively-oriented orthonormal basis
for $T_{p}\hat{M}$. Given an orientation of $\hat{M}$ this in turn specifies a canonical pseudohermitian 
structure $j:\xi^{\perp}\rightarrow\xi^{\perp}$, a feature which arises uniquely in dimension three. 

Now consider the product $\pi:N := \hat{M}\times(0,\infty)\rightarrow\hat{M}$, endowed with the structure 
of a {\em Riemannian cone} via the extended metric 
\[\bar{g} = t^{2}g + \frac{1}{4f} dt\otimes dt \hspace{.1in} \hbox{for}\hspace{.1in} t\in(0,\infty) \ , \ f := Ric(\xi) \ .\]
In a neighbourhood $U$ of any point $p\in\hat{M}$ we may choose a unit vector field $\sigma\in C^{\infty}(U,
\xi^{\perp})$ and define an orthonormal framing $\{\frac{1}{t}\sigma, \frac{1}{t}j\sigma, \frac{1}{t}\xi, 
2\sqrt{f}\frac{\partial}{\partial t}\}$ for $T^{*}N$ over $\pi^{-1}(U)$. At the same time, we may extend 
$\pi^{*}j$ to an endomorphism $J:TN\rightarrow TN$ such that 
\[ J(\frac{1}{t}\xi) = 2\sqrt{f}\frac{\partial}{\partial t} \ ; J(\frac{\partial}{\partial t}) = 
\frac{-1}{2t\sqrt{f}}\xi  \ ,\]
with $J\mid_{\pi^{*}\xi^{\perp}} = \pi^{*}j$. $T^{1,0}N$ and $T^{0,1}N$ then correspond to the 
$\pm\bf{i}$-eigenspaces of the ${\Bbb C}$-linear extension of $J$, endowing $N$ with an almost-complex
structure in the usual way. Integrability (i.e., Frobenius-involutivity) of this structure is then well-known
 to be equivalent 
to the assumption ${\mathcal L}_{\xi}j = 0$ on $\hat{M}$, i.e., for any $\sigma\in C^{\infty}(U,\xi^{\perp})$,
the Lie bracket satisfies
\[ [\xi, j\sigma] = j [\xi, \sigma] \ .\]
In smoothly preserving $j$, we say that the Reeb flow is {\em conformal} \cite{HP}. Let the complex surface
so obtained from $N$ be denoted $X$, and let $\Omega$ denote the natural hermitian form on the holomorphic tangent bundle, $TX$, corresponding to ${\Bbb C}$-linear extension of the contraction of $\bar{g}$ with $J$. 
More specifically, 
\[ \Omega = \iota_{J} \bar{g} = t^{2}\pi^{*}\omega + \frac{t}{\sqrt{f}}\vartheta\wedge dt \qquad \ (1),\]
where $\omega\mid_{\xi^{\perp}}:=\iota_{j}g$, $\iota_{\xi}\omega := 0$. Using the basic properties of the
Levi-Civita connection on $T\hat{M}$, we note that for $u , v\in \xi^{\perp}$
\[ d\vartheta(u,v) = u(\vartheta(v)) - v(\vartheta(u)) - \vartheta([u,v]) \]
\[ \hspace{.2in} = u (g(\xi, v)) - v (g(\xi, u)) - g(\xi, [u,v]) \]
\[ \hspace{.2in} = g(D_{u}\xi, v) - g(D_{v}\xi, u) + g(\xi, Tor(u,v))\]
\[ \hspace{.2in} = (\iota_{\beta - \beta^{*}}g)(u,v) \ ,\]
where $\beta:\xi^{\perp}\rightarrow\xi^{\perp}$ is the linear operator defined by $\beta(v) = D_{v}\xi$.
We now recall from \cite{HP}, Lemma 3, the identity
\[-2Ric(\xi) = trace(\beta^{2}) + \xi(trace(\beta)) \ .\]
By \cite{HP}, Proposition 2, the integrability condition ${\mathcal L}_{\xi} j = 0$ is equivalent to 
$j\beta = \beta j$. If in addition it is assumed that the (geodesible) flow of $\xi$ is volume-preserving, i.e., 
$div(\xi) = trace(\beta) = 0$, then we have
\[ \beta = \sqrt{f}\cdot j \ , \hspace{.1in}\hbox{hence}\hspace{.1in} \beta - \beta^{*} = 2\sqrt{f}\cdot j \ ,\]
and so 
\[d\vartheta\mid_{\xi^{\perp}} = 2\sqrt{f}\cdot(\iota_{j}g) = 2\sqrt{f}\cdot(\omega\mid_{\xi^{\perp}}) \ .\] 
Since both of these forms vanish under contraction with $\xi$, we may make the identification
\[ d\vartheta = 2\sqrt{f}\cdot\omega \ \hspace{.1in} \hbox{on} \hspace{.1in} T\hat{M} \ .\] 

We remark that the assumption that the flow of $\xi$ is both conformal and 
volume-preserving is equivalent to $\xi$ being a Killing vector field, via the equation 
\[ {\mathcal L}_{\xi}g\mid_{\xi^{\perp}} = (div(\xi)) g \]
derived in \cite{HP}, Proposition 2. Finally, suppose there is a compact 
subset $K\subset\subset M$ outside which the function $f = Ric(\xi)$ is equal to a positive constant $c$. Returning to equation (1) above, we note that a simple comparison of types in this case implies that 
$d\Omega\mid_{\pi^{-1}(\hat{M}\setminus K)} = 0$ if and only if $2\sqrt{c}\cdot\omega = d\vartheta$. It follows 
that $\Omega$ corresponds to a K\"{a}hler form on $\pi^{-1}(\hat{M}\setminus K)\subseteq X$, i.e., 
$\hat{M}\setminus K$ is a Sasakian three-manifold (cf. e.g., \cite{BiHu}).

We summarise with the following

\begin{proposition} Let $\xi$ be a geodesible Killing unit vector field on the compact Riemannian 
three-manifold
$\hat{M}$, such that $Ric(\xi) > 0$ and is constant outside a compact subset $K$ of $M$. Then $\hat{M}
\setminus K$ is Sasakian. 
\end{proposition} 

We now introduce a smooth complex vector bundle $E\rightarrow M$, equipped with a connection
\[ \nabla : C^{\infty}(E)\rightarrow C^{\infty}(E\otimes T^{*}M) \ ,\]
and consider the pullback $\pi^{*}(E)\rightarrow N\setminus\pi^{-1}(P_{0})$. A volume form on
$N$, compatible with both the Riemannian cone-structure and the K\"{a}hler structure on 
$\pi^{-1}(\hat{M}\setminus K)\subseteq X$, is naturally defined by
$dVol_{\bar{g}} := \frac{1}{2}\Omega\wedge\Omega$, such that 
\[ dVol_{\bar{g}}\mid_{\pi^{-1}(\hat{M}\setminus K)} = \frac{t^3}{2c}d\vartheta\wedge\vartheta\wedge
dt \ ,\]
noting that the induced volume-form/orientation on $\hat{M}$, corresponding to the slice $\{t = 1\}$, 
will then be 
\[\iota_{2\sqrt{c}\frac{\partial}{\partial t}}dVol_{\bar{g}} = \frac{-1}{\sqrt{c}}d\vartheta\wedge
\vartheta \ .\]
This entails a relative sign-change between the action of the Hodge star-operator $\star '$ on
$\bigwedge^{2}T^{*}N\mid_{\pi^{-1}(\hat{M}\setminus K)}$ and the corresponding operator $\star$ 
on $\bigwedge^{2}T^{*}\hat{M}$. In particular, 
if $F_{\nabla}$ denotes the curvature 2-form associated with $(E, \nabla)$, then we have
\[ \star ' (\pi^{*}F_{\nabla}) = \frac{-1}{2t\sqrt{c}}(\star\pi^{*}F_{\nabla})\wedge dt \qquad (2) .\]   

Choose a globally
defined endomorphism $\varphi\in C^{\infty}(M, E\otimes E^{*})$, and with it a natural extension
of $\nabla$ acting on sections of the pullback $\pi^{*}(E)\mid_{\pi^{-1}(\hat{M}\setminus K)}$, i.e.,
\[ \nabla ' := d + \pi^{*}(A_{\nabla}) + \frac{1}{t\sqrt{c}}\pi^{*}(\varphi)\otimes dt \ ,\]
where $A_{\nabla}$ denotes the connection matrix of $\nabla$ with respect to an arbitrary smooth 
framing of $E$. Consequently, we have 
\[ F_{\nabla '} = \pi^{*}F_{\nabla} + \frac{1}{t\sqrt{c}}\nabla\varphi\wedge dt \ ,\]
so that the identification
\[ F_{\nabla '} = \pi^{*}F_{\nabla} + \frac{1}{2t\sqrt{c}}(\star F_{\nabla})\wedge dt = 
\pi^{*}F_{\nabla} - \star '(\pi^{*}F_{\nabla})\]
follows automatically from the condition
\[ \nabla\varphi = \frac{1}{2}\star F_{\nabla} \qquad (3) \ .\]
Equation (3) is of course the {\em monopole equation} for the ensemble $(E,\nabla, \varphi)$ over 
$M$ \cite{M}, from which we recover the well-known $\star '$-anti self-duality condition for the curvature 2-form 
$F_{\nabla '}$ relative to $(\pi^{*}(E), \nabla ')\mid_{\pi^{-1}(M\setminus K)}$, and hence
a holomorphic structure on the bundle $\pi^{*}(E)\mid_{\pi^{-1}(M\setminus K)}$ \cite{A1}.  

The central problem of this discussion is to provide a sufficient condition for smooth completion of
the monopole data $(E,\nabla, \varphi)$, defined initially on $M$, over the compact Riemannian manifold
$\hat{M}$ with unit Killing vector field $\xi$ as described in Proposition 1. Our approach will be to 
formulate the sufficiency in terms of existence of a {\em holomorphic} connection 
\[{\mathcal D} : {\mathcal O}(\pi^{*}(E)\mid_{\pi^{-1}(M\setminus K)})\rightarrow\Omega^{1}_{X}
\otimes{\mathcal O}(\pi^{*}(E)\mid_{\pi^{-1}(M\setminus K)}) \ .\]
A ``Hartogs method" for extension of holomorphic vector bundles across gap-loci of sufficiently large
codimension, via holomorphic parallel transport of frames, was introduced in \cite{BuHa}, and adapted
to completion of monopole fields at point singularities in ${\Bbb R}^{3}$ in \cite{H}. Existence
of a holomorphic connection is well-known to be obstructed by analytic cohomology with coefficients
in the sheaf $\Omega^{1}_{X}\otimes{\mathcal O}(End(\pi^{*}(E)\mid_{\pi^{-1}(M\setminus K)}))$ \cite{A2}. 
As we shall
see, for the purposes of parallel transport a {\em relative} holomorphic connection is all that is required,
for which the obstruction may be formulated in terms of a Cauchy-Riemann equation of the form 
\[\bar{\partial}\Psi = \iota_{Z}F_{\nabla'} \ ,\hspace{.1in} \hbox{for} \hspace{.1in} \Psi\in
C^{\infty}(\pi^{-1}(M\setminus K), End(\pi^{*}(E)\mid_{\pi^{-1}(M\setminus K)})) \ ,\]
and where $Z$ denotes a holomorphic vector field on $X$. In order to formulate the sufficiency 
expressly in terms of the three-dimensional structure $(E, \nabla, \varphi)$, we will assume moreover
that $\Psi = \pi^{*}\psi$ for some $\psi\in C^{\infty}(M\setminus K, E\otimes E^{*})$, and similarly
that $Z$ corresponds to the trivial lifting to $\pi^{-1}(M\setminus K)$ of a smooth section of ${\Bbb C}\otimes\xi^{\perp}$, i.e., we will define $Z := \sigma - {\bf i}J\sigma$. $Z$ consequently lies in $T^{1,0}
N$ at each point of its definition, but must also be annihilated by the Cauchy-Riemann operator of the 
holomorphic tangent bundle $TX$.
If $D'$ denotes the Levi-Civita connection associated with $\bar{g}$ on $TN$, we may identify this Cauchy-Riemann 
operator as $\bar{\partial} = D' + {\bf i}JD'$, noting the important general fact that because the complex structure is also K\"{a}hler, we have $D'J = 0$. From this it follows easily that $\bar{\partial}Z = 0$.
Moreover, under these assumptions the Cauchy-Riemann equation for the obstruction $\iota_{Z}F_{\nabla'}$
may be reduced to three dimensions as follows. A $\bar{\partial}$-operator for the 
holomorphic vector bundle corresponding to $\pi^{*}(E)\mid_{\pi^{-1}(M\setminus K)}$ may similarly be 
defined in relation to $\nabla '$ by  $\bar{\partial}_{E}:= \nabla' + {\bf i}J\nabla'$, so that we may write
\[ \bar{\partial}_{E}\pi^{*}\psi = \bar{\partial}_{E}^{\perp}\pi^{*}\psi + (\nabla_{\xi} - 2{\bf i}\pi^{*}\varphi)
\pi^{*}\psi\otimes\vartheta + \frac{1}{t\sqrt{c}}(\pi^{*}\varphi + \frac{\bf{i}}{2}\nabla_{\xi})\pi^{*}\psi\otimes dt \qquad (4) ,\]
where $\bar{\partial}_{E}^{\perp} = \nabla^{\perp} + {\bf i}J\nabla^{\perp}$ denotes the natural restriction
of the operator to ${\Bbb C}\otimes\pi^{*}\xi^{\perp}$. Now, by a simple rearrangement of terms,
\[ \bar{\partial}_{E}^{\perp}\pi^{*}\psi = \nabla_{\sigma + {\bf i} J\sigma}\psi\otimes(\sigma^{*}
 - {\bf i}(j\sigma)^{*}) = \nabla_{\bar{Z}}\psi\otimes(\sigma^{*} - {\bf i}(j\sigma)^{*}) \ ,\]
where, for convenience, we will continue to use ``$Z$" to denote $\sigma - {\bf i}J\sigma$ when considered 
simply as a section of ${\Bbb C}\otimes T\hat{M}$.  
Relative to the same local framing of $T^{*}\hat{M}$ we now write
\[F_{\nabla} = F_{1,2}\sigma^{*}\wedge(j\sigma)^{*} + F_{1,3}\sigma^{*}\wedge\vartheta + F_{2,3}(j\sigma)^{*}
\wedge\vartheta \ ,\]
so that 
\[ (\iota_{Z}F_{\nabla '})^{\perp} = (\iota_{Z}\pi^{*}F_{\nabla})^{\perp} = {\bf i}F_{1,2}(\sigma^{*} - {\bf i}(j\sigma)^{*}) \  \ \text{i.e.,}\]
\[\nabla_{\bar{Z}}\psi = {\bf i}F_{1,2} \ .\]
From the monopole equation we obtain, on the other hand, 
\[ \nabla\varphi = \frac{1}{2}(F_{1,2}\vartheta - F_{1,3}(j\sigma)^{*} + F_{2,3}\sigma^{*}) \ ,\]
hence $\nabla_{\xi}\varphi = \frac{1}{2}F_{1,2}$. In conclusion,
\[\bar{\partial}_{E}^{\perp}\pi^{*}\psi = (\iota_{Z}F_{\nabla '})^{\perp}\hspace{.1in}\text{if and only if}
\hspace{.1in} \nabla_{\bar{Z}}\psi = 2{\bf i}\nabla_{\xi}\varphi \ .\]
Turning now to the transversal terms of equation (4), we recall that
\[F_{\nabla '} = \pi^{*}F_{\nabla} + \frac{1}{t\sqrt{c}}(\nabla\varphi)\wedge dt \ ,\]
hence
\[ (\iota_{Z}F_{\nabla '})_{\vartheta} = (\iota_{Z}F_{\nabla})_{\vartheta} = (F_{1,3} - {\bf i}F_{2,3}) = -2{\bf i}\nabla_{Z}\varphi \ .\]
We consequently derive a pair of identifications
\[(\nabla_{\xi} - 2{\bf i}\varphi)\psi = -2{\bf i}\nabla_{Z}\varphi \hspace{.1in}\text{and}\]
\[\frac{1}{t\sqrt{c}}(\varphi + \frac{\bf i}{2}\nabla_{\xi})\psi = \frac{1}{t\sqrt{c}}\nabla_{Z}\varphi \ ,\]
which are clearly the same equation. In conclusion, 
\begin{proposition} The equation
\[\bar{\partial}_{E}\pi^{*}\psi = \iota_{Z}F_{\nabla '} \]
is equivalent to the three-dimensional coupled system
\[\nabla_{\bar{Z}}\psi = 2{\bf i}\nabla_{\xi}\varphi\qquad ;\qquad (\nabla_{\xi} - 2{\bf i}\varphi)\psi 
= -2{\bf i}\nabla_{Z}\varphi \ .\]
Letting $\phi := -2{\bf i}\varphi$ and $\hat{\nabla} := \nabla + \phi\otimes\vartheta$, we may write this more concisely in the form
\[ \hat{\nabla}_{\bar{Z}}\psi = -\hat{\nabla}_{\xi}\phi\]
\[ \hat{\nabla}_{\xi}\psi = \hat{\nabla}_{Z}\phi \ .\]
\end{proposition} 

We will take a closer look at this system in the next section.

\vspace{.1in} 

\section{Higgs Potential of the coupled system}

\vspace{.1in}

The tensor $\varphi$ which comes as part the monopole data $(E, \nabla, \varphi)$ on the three-manifold
$M$, commonly referred to as the {\em Higgs Field}, will be said to admit a {\em Higgs Potential}
relative to the flow of a vector field $\eta$ if the equation
\[ \hat{\nabla}_{\eta}\Phi = \varphi \]
is solvable on $M$. If all data are smooth, and $\eta$ is non-vanishing, then local existence of a 
solution $\Phi$ reduces to an elementary ODE problem. More specifically, given $p\in M$ and a sufficiently
small neighbourhood $U_{p}$ in which the flow of $\eta$ is everywhere transversal to a disc $\Delta(p)$
centred at $p$, we may first take a smooth framing of $E\mid_{\Delta(p)}$, which supplies the initial
data for parallel transport along the flow, and hence for a smooth framing of $E\mid_{U_{p}}$ in which 
our equation is uniquely solvable, after specifying, e.g., the condition
\[\Phi\mid_{\Delta(p)} \ \equiv \ 0 \ .\]
 Diffeomorphic transport of the disc $\Delta(p)$ along the flow of 
$\eta$ moreover allows smooth continuation of $\Phi$ as far as possible throughout a chain of such neighbourhoods. Now consider the unit Killing field $\xi$ of the previous section, and let us examine the problem when the data $(E, \nabla, \varphi)$ are restricted to $M\setminus K$, which we may effectively identify with a punctured neighbourhood of $P_{0}\in\hat{M}$, for some compact subset $K$ of $M$. Again,
we take $\Delta(P_{0})$ to be a disc transversal to the flow of $\xi$ inside the complete neighbourhood $U_{P_{0}}$ in $\hat{M}$, which we may assume to be diffeomorphic to a product $\Delta(P_{0})\times(-\varepsilon,
\varepsilon)$. It is then straightforward to derive the existence of a solution to the equation
\[\hat{\nabla}_{\xi}\Phi = \varphi \]
inside $U_{P_{0}}\setminus\overline{B_{\delta}(P_{0})}$, such that
\[\Phi\mid_{\Delta(P_{0})} \ \equiv \  0\]
for $B_{\delta}(P_{0})$ a geodesic ball of arbitrarily small radius $\delta > 0$ centred at $P_{0}$.
If we now define
\[\psi := -2{\bf i}\hat{\nabla}_{Z}\Phi \ ,\]
we have a solution of the coupled system of Proposition 2, in $U_{P_{0}}
\setminus\{P_{0}\} = M\setminus K$, provided
\[\hat{\nabla}_{\bar{Z}}(\hat{\nabla}_{Z}\Phi) + \hat{\nabla}_{\xi}(\hat{\nabla}_{\xi}\Phi) = 0 \ ,\]
and
\[\hat{\nabla}_{Z}(\hat{\nabla}_{\xi}\Phi) = \hat{\nabla}_{\xi}(\hat{\nabla}_{Z}\Phi) \ .\]    
Via a standard expansion of second-order covariant differentiation of endomorphisms, we note that
\[ \hat{\nabla}_{Z}(\hat{\nabla}_{\xi}\Phi) - \hat{\nabla}_{\xi}(\hat{\nabla}_{Z}\Phi) = 
\hat{\nabla}_{[Z,\xi]}\Phi + [\iota_{Z}\iota_{\xi}F_{\hat{\nabla}}, \Phi] \qquad\qquad (5) ,\]
where the bracket $[Z,\xi]$ denotes the ${\Bbb C}$-linear extension of the Lie bracket of vector
fields, while $[\iota_{Z}\iota_{\xi}F_{\hat{\nabla}}, \Phi]$ corresponds to the formal commutator 
of matrices. If we now recall that
\[\hat{\nabla} = \nabla + \phi\otimes\vartheta \]
it is a direct calculation to show that
\[\iota_{Z}\iota_{\xi}F_{\hat{\nabla}} = \iota_{Z}\iota_{\xi}F_{\nabla} - \nabla_{Z}\phi \ .\]
Now $F_{\nabla} = 2\star\nabla\varphi$ implies 
\[\iota_{Z}\iota_{\xi}F_{\nabla} = 2\iota_{Z}\iota_{\xi}((\nabla_{\sigma}\varphi)\otimes(j\sigma)^{*}
\wedge\vartheta - (\nabla_{j\sigma}\varphi)\otimes\sigma^{*}\wedge\vartheta + (\nabla_{\xi}\varphi)\otimes
\sigma^{*}\wedge(j\sigma)^{*})\]
\[ = -2(-{\bf i}\nabla_{\sigma}\varphi - \nabla_{j\sigma}\varphi) = 2{\bf i}\nabla_{Z}\varphi = -\nabla_{Z}
\phi  \ ,\]
and consequently $\iota_{Z}\iota_{\xi}F_{\hat{\nabla}} = -2\nabla_{Z}\phi$. 

Recalling now the natural properties of the Levi-Civita connection $D$ on $TM$, together with the 
geodesibility of the Reeb vector field $\xi$, and the assumption $D_{\sigma}\xi = \beta(\sigma) = 
\sqrt{c}\cdot j\sigma$ on $M\setminus K$, we see that
\[ g([\sigma,\xi], \xi) = -g(D_{\sigma}\xi, \xi) = -\sqrt{c}\cdot g(j\sigma, \xi) = 0 \ ,\]
and 
\[ g([\sigma, \xi], \sigma) = g(D_{\sigma}\xi, \sigma) - g(D_{\xi}\sigma, \sigma) \ ,\]
\[ \qquad \qquad = \sqrt{c}\cdot g(j\sigma, \sigma) + g(\sigma, D_{\xi}\sigma) \ ,\]
hence $g([\sigma, \xi], \sigma) = 0$. This leaves us to conclude that $[\sigma, \xi] = 
\lambda j\sigma $ for some ${\Bbb R}$-valued function $\lambda = g([\sigma, \xi], j\sigma) $, and, with the additional integrability condition $[\xi, j\sigma] = j[\xi, \sigma]$, 
that $ [Z, \xi] = {\bf i}\lambda Z $. The commutator of covariant derivatives (5) may thus be rewritten
\[ \hat{\nabla}_{Z}(\hat{\nabla}_{\xi}\Phi) - \hat{\nabla}_{\xi}(\hat{\nabla}_{Z}\Phi) = 
{\bf i}\lambda\nabla_{Z}\Phi + 4{\bf i}[\nabla_{Z}\varphi, \Phi] \ .\]

A sufficient condition for solvability of the coupled equation of Proposition 2 may now be stated as follows.

\begin{proposition} If the Higgs Potential $\Phi$ relative to the Reeb vector field $\xi$ above satisfies 
the equations
\[ \nabla_{\bar{Z}}(\nabla_{Z}\Phi) + \nabla_{\xi}\varphi = 0 \ ,\]
and
\[ \lambda\nabla_{Z}\Phi + [\nabla_{Z}\varphi, \Phi]  = 0\]
for a unit vector field $\sigma\in C^{\infty}(M\setminus K, \xi^{\perp})$ and a real-valued function $\lambda = \frac{1}{4}g([\sigma, \xi], j\sigma) $, then $\psi := -2{\bf i}\hat{\nabla}_{Z}\Phi$ satisfies the coupled system 
\[ \hat{\nabla}_{\bar{Z}}\psi = -\hat{\nabla}_{\xi}\phi\]
\[ \hat{\nabla}_{\xi}\psi = \hat{\nabla}_{Z}\phi \ .\]
\end{proposition} 

\vspace{.1in}

{\em Remark:} Following \cite{BiHu}, we note that for the special case of $(\hat{M}, g)$ a {\em regular} 
Sasaki manifold, in particular, a compact three-manifold on which the flow of the Reeb vector field $\xi$ induces
the structure of an ${\Bbb S}^{1}$-principal bundle over a compact Riemann surface, it may be assumed that 
the Lie bracket $[Z, \xi] = 0$, and hence $\lambda = 0$ in the second equation above.  
   
\vspace{.1in}

At this point we recall the motivation for Propositions 2 and 3, namely to solve the Cauchy-Riemann
equation $\bar{\partial}_{E}\pi^{*}\psi = \iota_{Z}F_{\nabla '}$ in $\pi^{-1}(M\setminus K)$ thereby
removing the obstruction to existence of a relative holomorphic connection on $\pi^{*}E\mid_{\pi^{-1}(M\setminus K)}$.
The role of holomorphic connections in defining unique extension of the holomorphic bundle $\pi^{*}E$
to $\pi^{-1}(\hat{M}\setminus K)$, and hence extension of $E$ to $\hat{M}$, will be recalled in the next section.

\vspace{.1in}

\section{A Hartogs technique for point-completion of $E$}

\vspace{.1in}

Let us begin by once again identifying $M\setminus K$ with a punctured neighbourhood $U\setminus\{P_{0}\}$
in $\hat{M}$. Note that in a neighbourhood $U'$ of $P_{0}\in\pi^{-1}(U)$ we are in a position to choose holomorphic coordinates $(z,w)$ such that the $t$-axis, corresponding to the real line $\pi^{-1}(P_{0})$, lies wholly in the complex line defined in $U'$ by $z = 0$, and of course contains the point $P_{0} = (0,0)$. 
Consider a sufficiently small bi-disc $\Delta_{\varepsilon}\times\Delta_{\varepsilon}$ centred at the origin of this local holomorphic coordinate neighbourhood, where $\varepsilon$ denotes the radius in each coordinate, and the domains 
\[V := (\Delta_{\varepsilon}\setminus\Delta_{\frac{\varepsilon}{2}})\times\Delta_{\varepsilon}  \ , \ \hspace{.1in}\text{and}\hspace{.1in} W := \Delta_{\varepsilon}\times D \ ,\]
where $D$ is an open disc contained in $\Delta_{\varepsilon}\setminus\pi^{-1}(P_{0})$. From these we construct a Hartogs figure $H := V\cup W$. Note that on each of the components of $H$ we may appeal to the Oka-Grauert Principle for existence of a local holomorphic framing of $\pi^{*}E$ - let us refer to these as ${\bf f}_{V}$ and ${\bf f}_{W}$ respectively. Moreover, on $V\cap W$ suppose we have the transition relation ${\bf f}_{W} = \tau\cdot{\bf f}_{V}$ with respect to a holomorphic gauge-transformation $\tau$. Relative to ${\bf f}_{V}$, we denote the connection matrix $A^{1,0}_{Z} := \pi^{*}(A_{\sigma} - {\bf i}A_{j\sigma})$, while the corresponding component $A^{0,1}_{Z} = \pi^{*}(A_{\sigma} + {\bf i}A_{j\sigma})$, associated with the Cauchy-Riemann operator on $\pi^{*}E$, can be taken to be zero. Consequently,
\[ \iota_{Z}F_{\nabla'} = \bar{\partial}_{E}A_{Z}^{1,0} \ ,\]
so that the matrix $A_{Z}^{1,0} - \pi^{*}\psi$ is holomorphic. If we now choose $z_{0}\in\Delta_{\varepsilon}\setminus\Delta_{\frac{\varepsilon}{2}}$ the restriction ${\bf f}_{V}\mid_{\{z_{0}\}\times\Delta_{\varepsilon}}$ will provide the ``initial data" for a first-order system of holomorphic differential equations
\[  \frac{\partial{\bf f}}{\partial z} + (A_{Z}^{1,0} - \pi^{*}\psi){\bf f} = {\bf 0} \qquad (6)\]
in one complex variable. The solution may be extended via uniqueness of analytic continuation around each annular fibre of the product corresponding to $V$, with the sole obstruction to global continuation governed by a holonomy matrix at each point of $\{z_{0}\}\times\Delta_{\varepsilon}$. A similar holomorphic system may be defined on $W$, relative to the initial data provided by ${\bf f}_{W}\mid_{\{z_{0}\}\times D}$ - an important difference being that a covariantly constant solution $\tilde{\bf f}$ is unobstructed, given that the fibres of the product corresponding to $W$ are simply connected. 

Let $\gamma: \Delta_{\varepsilon}\times [0,1]\rightarrow V$ denote a family of loops, holomorphically 
parametrized by $w$ and traversing each annular fibre of $V$, such that $\gamma(w, 0) = \gamma(w, 1) = (w,z_{0})$. If $\alpha(w)$ denotes the holonomy matrix associated with analytic continuation of solutions to (6), then we may write ${\bf f}(\gamma(w, 1)) = \alpha\cdot {\bf f}(\gamma(w,0))$ and hence, for $w\in D$,
\[ \tilde{\bf f}(\gamma(w, 1)) = \tau\cdot\alpha\cdot \tau^{-1}\tilde{\bf f}(\gamma(w,0)) \ .\]
Because the fibres of $W$ are simply connected, $\alpha\mid_{D} \equiv {\bf 1}$, but $\alpha$ is holomorphic in $w$, hence we conclude that the holonomy everywhere in $V$ is trivial. The covariantly constant holomorphic framing of $\pi^{*}E\mid_{V}$ so obtained may easily be extended to $U'\setminus\pi^{-1}(P_{0})$ using analytic continuation along homotopy equivalent paths, and hence we arrive at a unique holomorphic extension of $\pi^{*}E$ to $U'$ via the classical Hartogs extension of functions. It follows at once that we have 

\vspace{.1in}

\begin{proposition} If the equations of Proposition 3 are satisfied by $\varphi$ and $\Phi$ in $M\setminus K$, then the vector bundle $E$ admits a smooth extension across $P_{0}$. 
\end{proposition}
  
\vspace{.1in}

\section{Completion of the Higgs field and monopole field}

\vspace{.1in}

Once a sufficient condition has been provided for smooth extension of the vector bundle $E$ over $\hat{M}$, a completion of the monopole field tensor corresponding to $F_{\nabla}$, as well as that of the Higgs field $\varphi$, can be addressed in a relatively straightforward way. Recall that a bounded smooth function $f$ defined on a punctured domain $U\setminus\{P_{0}\}$ in ${\Bbb R}^{n}$, for $n\geq 2$, need not admit a continuous extension at $P_{0}$ (let alone a smooth one). The problem of {\em boundary regularity}  is removed of course if we are prepared to assume that $f$ is uniformly H\"{o}lder continuous (or uniformly $C^{k,\alpha}$, if continuous extension of derivatives up to order $k$ is required) on $U\setminus\{P_{0}\}$. If simple functions are replaced by sections of a vector bundle $E$ defined over a Riemannian manifold $M$, some care needs to be taken in defining analogous notions of fractional differentiability for local frames of the bundle. Let ${\bf f}$ be such a smooth frame, defined over a domain ${\mathcal U}\subset M$. For an arbitrarily chosen base point $P\in{\mathcal U}$ we note that there is a smooth map $\rho:{\mathcal U}\rightarrow{\Bbb G}{\Bbb L}_{rk(E)}({\Bbb C})$ such that
${\bf f}(p) = \rho(p)\cdot{\bf f}(P)$, where $\rho(P) = id$. This map will serve as a local {\em gauge} for ${\bf f}$ relative to $P$, and we will say that ${\bf f}$ is 

\vspace{.1in}

(a) {\em admissible} if $\inf_{\mathcal U}|\det(\rho)| > 0$, and

(b) {\em uniformly H\"{o}lder $\alpha$-continuous on} ${\mathcal U}$ if the associated gauge $\rho$ has this property, i.e., there exists a constant $C > 0$ such that for all $p,q\in\Omega$
\[ \|\rho(p) - \rho(q) \| \leq C\cdot d(p,q)^{\alpha}\qquad \text{for some} \ 0 < \alpha\leq 1 \ ,\]
where $d(p,q)$ denotes the Riemannian metric distance between any two points. 

\vspace{.1in}

Similarly, we will say that
${\bf f}$ is {\em uniformly} $C^{1,\alpha}$ on ${\mathcal U}\subset\subset{\mathcal U}'$ if, for every smooth unit vector field $\zeta\in C^{\infty}({\mathcal U}', S\hat{M})$, the matrix-valued function $d\rho(\zeta)$ is uniformly 
H\"{o}lder $\alpha$-continuous. Now let ${\bf f}'$ be another admissible
smooth frame of $E$ on $\Omega$, hence there exists a gauge transformation 
$\tau:{\mathcal U}\rightarrow{\Bbb G}{\Bbb L}_{rk(E)}({\Bbb C})$ such that the gauge representing 
${\bf f}'$, $\rho' = \tau^{-1}\cdot\rho\cdot\tau$. It is not difficult to see that a product of admissible $C^{1,\alpha}$-gauge transformations is again admissible and $C^{1,\alpha}$. Moreover, for ${\mathcal U}$ of finite metric diameter, a $C^{1,\alpha}$-image $\tau({\mathcal U})$ will be a relatively compact domain contained in ${\Bbb G}{\Bbb L}_{rk(E)}({\Bbb C})$, hence the smooth involution $\iota$ acting on the general linear group by inversion will be uniformly Lipschitz when restricted to this image, and $\iota(\tau) = \tau^{-1}$ will again be $C^{1,\alpha}$ on ${\mathcal U}$. Such $\tau$ consequently form a subgroup of $C^{1}$-gauge transformations on ${\mathcal U}$, and hence the orbit of an admissible $C^{1,\alpha}$-gauge $\rho$ under the action of this subgroup via conjugation preserves the fractional 
differentiability. Conversely, if ${\bf f}'$ is assumed to be represented by a $C^{1,\alpha}$-gauge $\rho'$,
such that $\inf_{\mathcal U}|\det(\rho')| > 0$, it is a similarly elementary consequence that the gauge transformation between these frames is also $C^{1,\alpha}$. As a result, we see that the orbit of a given $C^{1,\alpha}$-frame ${\bf f}$ is independent of the choice of base-point $P\in{\mathcal U}$, and we may then simply refer to {\em the orbit} of admissible $C^{1,\alpha}$-frames of $E$ on this domain. 

If $\varphi\in C^{1,\alpha}({\mathcal U}, E\otimes E^{*}) = C^{1,\alpha}({\mathcal U}, 
\mathfrak{G}\mathfrak{L}_{rk(E)}({\Bbb C}))$ denotes the Higgs field of a monopole on $M$ with respect to a given admissible $C^{1,\alpha}$-frame, we may conclude using the same arguments that the action corresponding to the adjoint representation 
\[ Ad : G \rightarrow{\Bbb G}{\Bbb L}_{rk(E)}(\mathfrak{G}) \]
(where $G := {\Bbb G}{\Bbb L}_{rk(E)}({\Bbb C})$ and $\mathfrak{G}$ is the Lie Algebra),
when restricted to admissible $C^{1,\alpha}$-gauge transformations on ${\mathcal U}$, induces an orbit of $C^{1,\alpha}$-representatives of $\varphi$. In a similar manner, if $A$ denotes the connection matrix locally representing $\nabla$ in the same local framing, we will say that $\nabla$ is uniformly $C^{1,\alpha}$ on ${\mathcal U}$, and write $A\in C^{1,\alpha}({\mathcal U}, T^{*}M\otimes\mathfrak{G})$, if the contractions $dA(\zeta, \eta)$ are uniformly H\"{o}lder $\alpha$-continuous for every pair of smooth unit vector fields on ${\mathcal U}'$. This is moreover a property that is preserved under $C^{1,\alpha}$-gauge transformations. With respect to the derived adjoint representation 
\[ ad: \mathfrak{G} \rightarrow\mathfrak{G}\mathfrak{L}(\mathfrak{G}) \ ;\qquad\qquad \psi\mapsto [\psi, * ] \ ,\]
we note in particular that for every pair of smooth unit vector fields $\zeta \ , \ \eta$, the curvature 
\[F_{\nabla}(\zeta , \eta) = dA(\zeta , \eta) + [A_{\zeta}, A_{\eta}] = dA(\zeta , \eta) + 
ad(A_{\zeta})(A_{\eta}) \ .\]
Given $\nabla$ uniformly $C^{1,\alpha}$ on ${\mathcal U}\subset\subset\hat{M}$, $A_{\zeta}$ and $A_{\eta}$ are uniformly Lipschitz, hence it is a simple consequence of the triangle inequality, together with the skew-symmetry
\[ ad(A_{\zeta})(A_{\eta}) = -ad(A_{\eta})(A_{\zeta}) \ ,\]
that
\[ \| ad(A_{\zeta}(p))(A_{\eta}(p)) - ad(A_{\eta}(q))(A_{\zeta}(q)) \|\qquad\qquad\qquad\qquad \]
\[ \qquad\qquad \leq \|ad(A_{\zeta}(p))\| \|A_{\eta}(p) - A_{\eta}(q) \| + \| ad(A_{\eta}(q))\| \| A_{\zeta}(p) - A_{\zeta}(q) \| \ .\]

Since the linear operator norms of $ad(A_{\zeta})$ and $ad(A_{\eta})$ are themselves uniformly bounded, we conclude that the commutators $[A_{\zeta}, A_{\eta}]$ are uniformly Lipschitz on ${\mathcal U}$.

When ${\mathcal U}$ corresponds to a punctured neighbourhood $U\setminus\{P_{0}\}\subseteq\hat{M}$, supporting a frame of $E$ which extends smoothly to all of $U$, it follows directly from the assumption of uniform $C^{1,\alpha}$-regularity for $\varphi$ 
and $\nabla$ that both admit $C^{1}$-extensions across $P_{0}$ (and hence to $\hat{M}$). Moreover, the combination of uniform H\"{o}lder $\alpha$-continuity and Lipschitz continuity which these assumptions imply for the curvature, $F_{\nabla}$, allows continuous extension of this form across $P_{0}$ as well.

We may now summarize the results of this and the preceding sections in our main conclusion.

\begin{theorem} Let $(\hat{M}, g)$ be a compact, oriented Riemannian three-manifold corresponding to the point-completion of a manifold $M$, and let $\xi$ denote a geodesible Killing unit vector field on $\hat{M}$ such that the Ricci curvature function $Ric_{g}(\xi) > 0$ everywhere, and is constant outside a compact subset $K\subset M$. Suppose further that $(E,\nabla, \varphi)$ supply the essential data of a monopole field on $M$, smooth outside isolated singularities all contained in $K$, and that the following hold in the complement of $K$:

\vspace{.1in}

(i) the associated Higgs potential $\Phi$ satisfies the equations
\[ \nabla_{\bar{Z}}(\nabla_{Z}\Phi) + \nabla_{\xi}\varphi = 0 \ ,\]
and
\[ \lambda\nabla_{Z}\Phi + [\nabla_{Z}\varphi, \Phi]  = 0\]
where $Z = \sigma - {\bf i}j\sigma$ with respect to an orthonormal framing $\{\sigma, j\sigma, \xi\}\in C^{\infty}(\hat{M}\setminus K, T\hat{M})$, and a real-valued function $\lambda = \frac{1}{4}g([\sigma, \xi], j\sigma) $,

\vspace{.1in}

(ii) the Higgs field $\varphi$ and monopole connection $\nabla$ are both uniformly $C^{1,\alpha}$ with respect to admissible $C^{1,\alpha}$-framings of $E$.

\vspace{.1in}

Then the bundle $E$ extends smoothly across $P_{0}$, and admits $C^{1}$-extensions of both $\varphi$ and $\nabla$, together with a continuous extension of $F_{\nabla}$.  
 
\end{theorem} 

\vspace{.1in}

\begin{corollary} When $\hat{M} = {\Bbb S}^{3}$, equipped with the round metric, and the flow of $\xi$ induces the Hopf fibration, then the simplified condition

\[(i)\qquad  \nabla_{\bar{Z}}(\nabla_{Z}\Phi) + \nabla_{\xi}\varphi
 = [\nabla_{Z}\varphi, \Phi]  = 0 \ ,\]
  
together with condition $(ii)$, is sufficient for extension of $E, \varphi, \nabla$ and $F_{\nabla}$
as above.

\end{corollary}

\end{document}